\documentclass{article}
\usepackage{amssymb}
\usepackage{amsmath}

\newtheorem{theorem}{Theorem}

\newtheorem{corollary}[theorem]{Corollary}

\newtheorem{definition}[theorem]{Definition}

\newenvironment{proof}[1][Proof]{\noindent\textbf{#1.} }{\ \rule{0.5em}{0.5em}}
\input{tcilatex}
\begin{document}

\title{Slant helices and Darboux helices in Myller Configuration }
\author{Ak{\i}n Alkan$^1$, Mehmet \"{O}nder$^2$ \\
%EndAName
$^1$Manisa Celal Bayar University, \\
G\"{o}rdes Vocational School, 45750, G\"{o}rdes, Manisa, Turkey.\\
$^2$Delibekirli Village, Tepe Street, 31440, K{\i}r{\i}khan, Hatay, Turkey.\\
E-mails: $^1$akin.alkan@cbu.edu.tr, $^2$mehmetonder197999@gmail.com \\
Orcid Ids: $^1$https://orcid.org/0000-0002-8179-9525, \\
$^2$https://orcid.org/0000-0002-9354-5530}
\maketitle

\begin{abstract}
In this paper, we study slant helix (or $\overset{\_}{\xi}_{2}$-helix) and
Darboux helix in Myller configuration $M$. We show that a curve in $M$ is a
slant helix if and only if it is a Darboux helix. We give the alternative
frame of a curve in $M$. Furthermore, we obtain the differential equations
characterizing the curves in $M$ by means of both Frenet type frame and
alternative frame.
\end{abstract}

\textbf{AMS Classsification:} 53A04.


\section{\protect\bigskip Introduction}

\bigskip Orthonormal frames are the most useful tools to study the
differential geometry of curves and surfaces. Frenet-Serret frame of a space
curve and Darboux frame of a surface curve are the most famous ones of such
frames. More general frames along space curves or surface curves have been
introduced by Myller. By taking into account a unit vector field $\overset{\_%
}{\xi }$ and a plane $\pi $ along a curve $C$ and calling them as a versor
field $(C,\overset{\_}{\xi })$ and a plane field $(C,\pi )$, respectively,
such that~$\overset{\_}{\xi }\in \pi ,$ he defined a configuration in$~E^{3}$
called with his name as Myller configuration and denoted by $M(C,\overset{\_}%
{\xi },\pi )$\cite{Miron}. When the planes $\pi ~$are tangent to $C$, a
special case of this configuration is obtained and called a tangent Myller
configuration, denoted by $M_{t}(C,\overset{\_}{\xi },\pi )$. Especially, if
the curve $C$ is a surface curve lying on a surface $S\subset E^{3}$ with
arclength parameter $s$, $\overset{\_}{\xi }(s)$ is the tangent versor field
to $S$ along $C$, $\pi (s)$ is the tangent plane field to $S$ along $C$,
then $M_{t}(C,\overset{\_}{\xi },\pi )$ is called the tangent Myller
configuration intrinsic associated to the geometric objects $S,~C,$ $\overset%
{\_}{\xi }$. Then, the geometry of the versor field $(C,\overset{\_}{\xi })$
on a surface $S$ is the geometry of the associated Myller configurations $%
M_{t}(C,\overset{\_}{\xi },\pi )$ and moreover, $M_{t}(C,\overset{\_}{\xi }%
,\pi )$ represents a particular case of $M(C,\overset{\_}{\xi },\pi )$. In
the special case that tangent Myller configuration $M_{t}(C,\overset{\_}{\xi 
},\pi )$ is the associated Myller configuration to a curve $C$ on a surface $%
S,$ the classical theory of surface curves (curves lying on a surface) is
obtained.

\bigskip Myller studied the parallelism of versor field $(C,\overset{\_}{\xi 
})$ in the plane field $(C,\pi )$ \cite{Myller}. He obtained a
generalization of parallelism of Levi-Civita on the curved surfaces. Later,
Mayer gave some new invariants for $M(C,\overset{\_}{\xi },\pi )$ and $%
M_{t}(C,\overset{\_}{\xi },\pi )$ \cite{Mayer}. Miron extended the notion of
Myller configuration in Riemannian Geometry \cite{Miron}. Vaisman considered
Myller configuration in the symplectic geometry \cite{Vaisman}. Furthermore,
a Myller configuration was studied in differenet spaces \cite%
{Miron,Constantinescu}. Macsim, Mihai and Olteanu have studied
rectifying-type curves in a Myller configuration \cite{Macsim}. Recently,
the authors have defined some special helices in $M$ and given the
characterizations for those curves \cite{Akin}.

In the present paper, we study slant helices and Darboux helices in Myller
configuration $M$. We give characterizations and axes for these curves.
Moreover, we introduce the alternative frame of a curve $C$ in $M$ and give
the geometric interpretations of alternative invariants of $C$. Finally, we
obtain the differential equations related to curves, helices and slant
helices in $M$.

\section{Preliminaries}

In this section, a brief summary of curves in Myller configuration $M(C,%
\overset{\_}{\xi },\pi )$ are introduced. For more detailed information, we
refer to \cite{Miron}.

\bigskip Let $(C,\overset{\_}{\xi })$ be a versor field and $(C,\pi )$ be a
plane field such that~$\overset{\_}{\xi }\in \pi $ in $E^{3}$. The pair $((C,%
\overset{\_}{\xi }),(C,\pi ))$ is called Myller configuration in $E^{3}$ and
denoted by $M(C,\overset{\_}{\xi },\pi )$ \ or briefly $M$. Let $R=(O;%
\overset{\_}{i}_{1},\overset{\_}{i}_{2},\overset{\_}{i}_{3})$ denotes an
orthonormal frame. Then, $(C,\overset{\_}{\xi })$ can be written in the form

\begin{eqnarray*}
\overset{\_}{r} &=&\overset{\_}{r}(s),~\overset{\_}{\xi }=\overset{\_}{\xi }%
(s),~s\in I=(s_{1},s_{2}), \\
\overset{\_}{r}(s) &=&x(s)\overset{\_}{i}_{1}+y(s)\overset{\_}{i}_{2}+z(s)%
\overset{\_}{i}_{3}=\overrightarrow{OP}(s), \\
\overset{\_}{\xi }(s) &=&\overset{\_}{\xi }_{1}(s)\overset{\_}{i}_{1}+%
\overset{\_}{\xi }_{2}(s)\overset{\_}{i}_{2}+\overset{\_}{\xi }_{3}(s)%
\overset{\_}{i}_{3}=\overrightarrow{PQ},
\end{eqnarray*}%
where $\overset{\_}{r}=\overset{\_}{r}(s)$ is the position vector of $C$, $s$
is the arclength parameter of $C,$ and $\left\Vert \overset{\_}{\xi }%
(s)\right\Vert ^{2}=\left\langle \overset{\_}{\xi }(s),\overset{\_}{\xi }%
(s)\right\rangle =1.$ Writing $\overset{\_}{\xi }_{1}(s)=\overset{\_}{\xi }%
(s)$ and considering $\frac{d\overset{\_}{\xi }_{1}(s)}{ds},$ one can define
versor field $\overset{\_}{\xi }_{2}(s)$ as follows,%
\begin{equation*}
\frac{d\overset{\_}{\xi }_{1}(s)}{ds}=K_{1}(s)\overset{\_}{\xi }_{2}(s),
\end{equation*}%
where $K_{1}(s)=\left\Vert \frac{d\overset{\_}{\xi }}{ds}\right\Vert >0$ is
called curvature (or $K_{1}$-curvature) of $(C,\overset{\_}{\xi })$. The
versor field $\overset{\_}{\xi }_{2}(s)$ is orthogonal to $\overset{\_}{\xi }%
_{1}(s)$ and exists when $K_{1}(s)\neq 0$. By defining the versor field $%
\overset{\_}{\xi }_{3}(s)=\overset{\_}{\xi }_{1}(s)\times \overset{\_}{\xi }%
_{2}(s),$ an orthonormal and positively oriented frame is obained. This
frame is called the invariant frame of Frenet-type of the versor field $(C,%
\overset{\_}{\xi })$ and denoted by $R_{F}\left( P(s);\overset{\_}{\xi }%
_{1}(s),\overset{\_}{\xi }_{2}(s),\overset{\_}{\xi }_{3}(s)\right) $, or
briefly $R_{F}$ \cite{Miron}.

The derivative formulas for $R_{F}$ are

\begin{equation}
\frac{d\overset{\_}{r}(s)}{ds}=\overset{\_}{\alpha }(s)=a_{1}(s)\overset{\_}{%
\xi }_{1}(s)+a_{2}(s)\overset{\_}{\xi }_{2}(s)+a_{3}(s)\overset{\_}{\xi }%
_{3}(s),  \label{2.1}
\end{equation}%
with $a_{1}^{2}(s)+a_{2}^{2}(s)+a_{3}^{2}(s)=1$ and

\begin{eqnarray*}
\frac{d\overset{\_}{\xi }_{1}(s)}{ds} &=&K_{1}(s)\overset{\_}{\xi }_{2}(s),
\\
\frac{d\overset{\_}{\xi }_{2}(s)}{ds} &=&-K_{1}(s)\overset{\_}{\xi }%
_{1}(s)+K_{2}(s)\overset{\_}{\xi }_{3}(s), \\
\frac{d\overset{\_}{\xi }_{3}(s)}{ds} &=&-K_{2}(s)\overset{\_}{\xi }_{2}(s),
\end{eqnarray*}%
where $K_{2}(s)=\left\langle \frac{d\overset{\_}{\xi }_{3}(s)}{ds},\overset{%
\_}{\xi }_{3}(s)\right\rangle $ is called torsion (or $K_{2}$-torsion) of $%
(C,\overset{\_}{\xi }).$ The curvatures $K_{1}(s)$ and $K_{2}(s)$ have the
same geometrical interpretation as the same as the curvature and torsion of
a curve in $E^{3}$ and the functions $%
a_{1}(s),~a_{2}(s),~a_{3}(s),~K_{1}(s),~K_{2}(s),~(s\in I)~$are invariants
of the versor field $(C,\overset{\_}{\xi })$ \cite{Miron}. If $%
a_{1}(s)=1,~a_{2}(s)=0,~a_{3}(s)=0$, one obtains the Frenet equations of a
curve in $E^{3}.$

The following theorem is the fundamental theorem for the versor field $%
(C,\xi )$:

\begin{theorem}
(\cite{Miron})If the functions $%
K_{1}(s),K_{2}(s),a_{1}(s),a_{2}(s),a_{3}(s), $ $\left(
a_{1}^{2}(s)+a_{2}^{2}(s)+a_{3}^{2}(s)=1\right) $ of class $C^{\infty }$ are
a priori given, $s\in \lbrack a,b]$, then there exist a curve $%
C:[a,b]\rightarrow E^{3}$ with arclength $s$ and a versor field $\overset{\_}%
{\xi }(s)$, $s\in \lbrack a,b]$ such that the functions $a_{i}(s) $, $%
(i=1,2,3),$ $K_{1}(s)$ and$~K_{2}(s)$ are the invariants of $(C,\overset{\_}{%
\xi })$. Any two such versor fields $(C,\overset{\_}{\xi })$ differ by a
proper Euclidean motion.
\end{theorem}

We give the definition of $\overset{\_}{\xi }_{i}$-helix $(i=1,2,3)$ as
follows:\bigskip

\begin{definition}
(\cite{Akin}) A curve $C$ in $M$ with invarinat type Frenet frame $%
R_{F}\left( P(s);\overset{\_}{\xi }_{1}(s),\overset{\_}{\xi }_{2}(s),\overset%
{\_}{\xi }_{3}(s)\right) $ is called $\overset{\_}{\xi }_{i}$-helix if the
versor field $\overset{\_}{\xi }_{i}$ makes a constant angle with a constant
direction, where $i=1,2,3$.
\end{definition}

\begin{theorem}
(\cite{Akin}) A curve $C$ in $M$ is $\overset{\_}{\xi }_{1}$-helix (or $%
\overset{\_}{\xi }_{3}$-helix) according to the Frenet type frame $R_{F}$
iff \ $\frac{K_{2}}{K_{1}}$ is constant.
\end{theorem}

\section{Slant helices in Myller Configuration $M(C,\protect\overset{\_}{%
\protect\xi },\protect\pi )$}

In this section, we consider slant helix (or $\xi _{2}$-helix) in Myller
configuration $M(C,\overset{\_}{\xi },\pi )$ and we give the
characterizations for slant helix in $M$.

\begin{definition}
Let $(C,\overset{\_}{\xi })$ be a versor field with invariant type Frenet
frame $R_{F}\left( P(s);\overset{\_}{\xi }_{1}(s),\overset{\_}{\xi }_{2}(s),%
\overset{\_}{\xi }_{3}(s)\right) $. The curve $C$ is called a slant helix(or 
$\overset{\_}{\xi }_{2}$-helix) in $M$ if the versor field $\overset{\_}{\xi 
}_{2}$ makes a constant angle with a fixed unit direction $\overset{\_}{d%
\text{.}}$
\end{definition}

\begin{theorem}
\bigskip The curve $C$ in $M$ with Frenet frame $R_{F}$ and $\left(
K_{1},K_{2}\right) \neq (0,0)$ is slant helix iff the following function is
constant 
\begin{equation}
\sigma =\cot \theta =\mp \frac{K_{1}^{2}}{\left( K_{1}^{2}+K_{2}^{2}\right)
^{\frac{3}{2}}}\left( \frac{K_{2}}{K_{1}}\right) ^{\prime }.  \label{eqn3}
\end{equation}
\end{theorem}

\begin{proof}
From Definition 4, there exists a constant angle $\theta $ such that $%
\left\langle \overset{\_}{\xi }_{2},\overset{\_}{d}\right\rangle =\cos
\theta .$ Then, for the axis $\overset{\_}{d}$ we can write 
\begin{equation}
\overset{\_}{d}=x_{1}\overset{\_}{\xi }_{1}+(\cos \theta )\overset{\_}{\xi }%
_{2}+x_{3}\overset{\_}{\xi }_{3},  \label{3.1}
\end{equation}%
where $x_{1}=x_{1}(s),~x_{3}=x_{3}(s)$ are smooth functions of $s$. Since $%
\overset{\_}{d}$ is constant, differentiating (\ref{3.1}) with respect to $s$
gives 
\begin{equation*}
0=\left( x_{1}^{\prime }-K_{1}\cos \theta \right) \overset{\_}{\xi }%
_{1}+\left( x_{1}K_{1}-x_{3}K_{2}\right) \overset{\_}{\xi }_{2}+\left(
K_{2}\cos \theta +x_{3}^{\prime }\right) \overset{\_}{\xi }_{3},
\end{equation*}%
and, we have the system 
\begin{equation}
\left\{ 
\begin{array}{c}
x_{1}^{\prime }-K_{1}\cos \theta =0, \\ 
x_{1}K_{1}-x_{3}K_{2}=0, \\ 
K_{2}\cos \theta +x_{3}^{\prime }=0.%
\end{array}%
\right.   \label{3.2}
\end{equation}%
From the second equation in (\ref{3.2}), it follows $x_{1}=x_{3}\frac{K_{2}}{%
K_{1}}$. Writing that in the first and third equations in (\ref{3.2}) gives
differential equation 
\begin{equation}
x_{3}^{\prime }\left( 1+\left( \frac{K_{2}}{K_{1}}\right) ^{2}\right) +x_{3}%
\frac{K_{2}}{K_{1}}\left( \frac{K_{2}}{K_{1}}\right) ^{\prime }=0.
\label{3.3}
\end{equation}%
Considering the transformation $t=\frac{K_{2}}{K_{1}},$ the differential
equation in (\ref{3.3}) becomes $\frac{dx_{3}}{dt}\left( 1+t^{2}\right)
+x_{3}tt^{\prime }=0$ which has the solution $x_{3}=\lambda \frac{K_{1}}{%
\sqrt{K_{1}^{2}+K_{2}^{2}}}$, where $\lambda $ is integration constant.
Hence, $x_{1}=\lambda \frac{K_{2}}{\sqrt{K_{1}^{2}+K_{2}^{2}}}$. Since $%
\left\Vert \overset{\_}{d}\right\Vert =1$, from (\ref{3.1}) we have $\lambda
=\mp \sin \theta .$ Then, (\ref{3.1}) becomes%
\begin{equation}
\overset{\_}{d}=\mp \sin \theta \frac{K_{2}}{\sqrt{K_{1}^{2}+K_{2}^{2}}}%
\overset{\_}{\xi }_{1}+(\cos \theta )\overset{\_}{\xi }_{2}\mp \sin \theta 
\frac{K_{1}}{\sqrt{K_{1}^{2}+K_{2}^{2}}}\overset{\_}{\xi }_{3}.
\label{3.4.4}
\end{equation}%
From the third equation in (\ref{3.2}), it follows $\left( \mp \sin \theta 
\frac{K_{1}}{\sqrt{K_{1}^{2}+K_{2}^{2}}}\right) ^{\prime }=-\cos \theta
K_{2}.$ Hence we have that $\sigma =\cot \theta =\mp \frac{K_{1}^{2}}{\left(
K_{1}^{2}+K_{2}^{2}\right) ^{\frac{3}{2}}}\left( \frac{K_{2}}{K_{1}}\right)
^{\prime }$ is a constant function.

Conversely, let the function $\sigma $ given in (\ref{eqn3}) be constant and 
$\overset{\_}{d}$ be a unit vector defined by%
\begin{equation*}
\overset{\_}{d}=\mp \sin \theta \frac{K_{2}}{\sqrt{K_{1}^{2}+K_{2}^{2}}}%
\overset{\_}{\xi }_{1}+(\cos \theta )\overset{\_}{\xi }_{2}\mp \sin \theta 
\frac{K_{1}}{\sqrt{K_{1}^{2}+K_{2}^{2}}}\overset{\_}{\xi }_{3},
\end{equation*}%
where $\theta $ is constant. Differentiating last equalty gives 
\begin{equation*}
\overset{\_}{d}^{\prime }=K_{1}(\mp (\sin \theta )\sigma -\cos \theta )%
\overset{\_}{\xi }_{1}+K_{2}(\pm (\sin \theta )\sigma +\cos \theta )\overset{%
\_}{\xi }_{3}.
\end{equation*}%
Now, writing (\ref{eqn3}) in this result, we have $\overset{\_}{d}^{\prime
}=0$, i.e., $\overset{\_}{d}~$ is a constant vector field and since $%
\left\langle \overset{\_}{\xi _{2}},\overset{\_}{d}\right\rangle $ is
constant, we obtain that $C$ is a slant helix in Myller configuration $M.$
\end{proof}

From Theorem 5, the following corollary is obtained:

\begin{corollary}
\bigskip The axis of slant helix $C$ in Myller configuration $M$ is given by 
\begin{equation*}
\overset{\_}{d}=\mp \sin \theta \frac{K_{2}}{\sqrt{K_{1}^{2}+K_{2}^{2}}}%
\overset{\_}{\xi }_{1}+(\cos \theta )\overset{\_}{\xi }_{2}\mp \sin \theta 
\frac{K_{1}}{\sqrt{K_{1}^{2}+K_{2}^{2}}}\overset{\_}{\xi }_{3},
\end{equation*}%
where $\theta $ is the constant angle defined by $\left\langle \overset{\_}{%
\xi _{2}},\overset{\_}{d}\right\rangle =\cos \theta .$
\end{corollary}

\section{Darboux helices in Myller Configuration $M(C,\protect\overset{\_}{%
\protect\xi },\protect\pi )$}

Alternative frame of a curve is a useful tool to study the properties of the
curve in the Euclidean 3-space. This frame is obtained from the Frenet frame
of a curve and provide an advantage to give some characterizations of some
special curves and surfaces \cite{Kaya1,Kaya2,Sahiner,Uzunoglu}. In this
section, first we define alternative frame of a curve $C$ and give the
geometric interpretations of alternative curvatures. Next, we define Darboux
helices in $M$ and give characterizations of Darboux helices by means of
alternative frame.

Let $(C,\overset{\_}{\xi })$ be a versor field with invariant type Frenet
frame $R_{F}\left( P(s);\overset{\_}{\xi }_{1}(s),\overset{\_}{\xi }_{2}(s),%
\overset{\_}{\xi }_{3}(s)\right) $ in Myller configuration $M$. The vector $%
D=K_{2}\overset{\_}{\xi }_{1}+K_{1}\overset{\_}{\xi }_{3}$ is called Darboux
vector field of the Frenet frame $R_{F}$.

Let define $\overset{\_}{Y}=\frac{\overset{\_}{\xi }_{2}^{\prime }}{%
\left\Vert \overset{\_}{\xi }_{2}^{\prime }\right\Vert }$ and $\overset{\_}{D%
}=\frac{K_{2}}{\sqrt{K_{1}^{2}+K_{2}^{2}}}\overset{\_}{\xi }_{1}+\frac{K_{1}%
}{\sqrt{K_{1}^{2}+K_{2}^{2}}}\overset{\_}{\xi }_{3}$. Then, $R_{A}\{P(s);%
\overset{\_}{\xi }_{2}(s),\overset{\_}{Y}(s),\overset{\_}{D}(s)=\overset{\_}{%
\xi }_{2}\times \overset{\_}{Y}\}$ is an othonormal frame of the versor
field $(C,\overset{\_}{\xi })$ in $M$ \ and called alternative Frenet frame
(or AF-frame) of $(C,\overset{\_}{\xi })$. The moving equations of the
AF-frame are as follows%
\begin{equation}
\frac{d\overset{\_}{r}(s)}{ds}=\overset{\_}{\alpha }(s)=d_{1}(s)\overset{\_}{%
\xi }_{2}(s)+d_{2}(s)\overset{\_}{Y}(s)+d_{3}(s)\overset{\_}{D}(s),
\label{3.5}
\end{equation}%
with $d_{1}^{2}(s)+d_{2}^{2}(s)+d_{3}^{2}(s)=1$ and%
\begin{equation}
\begin{array}{c}
\overset{\_}{\xi }_{2}^{\prime }=p\overset{\_}{Y}, \\ 
\overset{\_}{Y}^{\prime }=-p\overset{\_}{\xi }_{2}+q\overset{\_}{D}, \\ 
\overset{\_}{D}^{\prime }=-q\overset{\_}{Y},%
\end{array}
\label{3.6}
\end{equation}%
where $d_{i}(s),~(i=1,2,3),~p=\sqrt{K_{1}^{2}+K_{2}^{2}}$ and $q=\frac{%
K_{1}^{2}}{K_{1}^{2}+K_{2}^{2}}\left( \frac{K_{2}}{K_{1}}\right) ^{\prime }$
are invariants and called alternative curvatures of the versor field $(C,%
\overset{\_}{\xi })$.

The geometric interpretations of the functions $d_{i}(s),~(i=1,2,3)$ can
simply be given as 
\begin{equation*}
d_{1}(s)=\cos \sphericalangle (\overset{\_}{\alpha },\overset{\_}{\xi }%
_{2}),~d_{2}(s)=\cos \sphericalangle (\overset{\_}{\alpha },\overset{\_}{Y}%
),~d_{3}(s)=\cos \sphericalangle (\overset{\_}{\alpha },\overset{\_}{D}).
\end{equation*}

For the geometric interpretations of the invariants $p,q,$ we give the
following theorem.

\begin{theorem}
Let consider a variation of alternative frame $R_{A}\{P(s);\overset{\_}{\xi }%
_{2}(s),\overset{\_}{Y}(s),\overset{\_}{D}(s)\}\rightarrow R_{A}^{\ast
}\{P^{\ast }(s+\Delta s);\overset{\_}{\xi }_{2}(s+\Delta s),\overset{\_}{Y}%
(s+\Delta s),\overset{\_}{D}(s+\Delta s)\}$ and let $\varphi _{1}$ be
oriented angle between the successive versor fields $\overset{\_}{\xi }%
_{2}(s)$, and $\varphi _{2}$ be oriented angle between the successive versor
fields $\overset{\_}{D}(s)$. Alternative curvatures $p$ and $q$ of the
versor field $(C,\overset{\_}{\xi })$ are given by%
\begin{equation}
p=\left\vert \frac{d\varphi _{1}}{ds}\right\vert ,~q=\left\vert \frac{%
d\varphi _{2}}{ds}\right\vert ,  \label{3.6*}
\end{equation}%
respectively.
\end{theorem}

\begin{proof}
Let $\varphi _{1}(s)$ be oriented angle between the successive versor fields 
$\overset{\_}{\xi }_{2}(s)$. Then, the angle function between versor fields $%
\overset{\_}{\xi }_{2}(s)$ and $\overset{\_}{\xi }_{2}(s+\Delta s)$ is
defined by $\left\vert \varphi _{1}(s+\Delta s)-\varphi _{1}(s)\right\vert
=h(s)$. From (\ref{3.6}), we have $p(s)=\left\Vert \overset{\_}{\xi }%
_{2}^{\prime }(s)\right\Vert $. Then, we can write 
\begin{eqnarray}
p(s) &=&\left\Vert \overset{\_}{\xi }_{2}^{\prime }(s)\right\Vert  \notag \\
&=&\lim_{\Delta s\rightarrow 0}\left\Vert \frac{\overset{\_}{\xi }%
_{2}(s+\Delta s)-\overset{\_}{\xi }_{2}(s)}{\Delta s}\right\Vert  \notag \\
&=&\lim_{\Delta s\rightarrow 0}\left\Vert \frac{\varphi _{1}(s+\Delta
s)-\varphi _{1}(s)}{\Delta s}\frac{\overset{\_}{\xi }_{2}(s+\Delta s)-%
\overset{\_}{\xi }_{2}(s)}{\varphi _{1}(s+\Delta s)-\varphi _{1}(s)}%
\right\Vert  \notag \\
&=&\lim_{\Delta s\rightarrow 0}\left\vert \frac{\varphi _{1}(s+\Delta
s)-\varphi _{1}(s)}{\Delta s}\right\vert \lim_{\Delta s\rightarrow
0}\left\Vert \frac{\overset{\_}{\xi }_{2}(s+\Delta s)-\overset{\_}{\xi }%
_{2}(s)}{\varphi _{1}(s+\Delta s)-\varphi _{1}(s)}\right\Vert  \notag \\
&=&\left\vert \frac{d\varphi _{1}}{ds}\right\vert \lim_{\Delta s\rightarrow
0}\left\Vert \frac{\overset{\_}{\xi }_{2}(s+\Delta s)-\overset{\_}{\xi }%
_{2}(s)}{\varphi _{1}(s+\Delta s)-\varphi _{1}(s)}\right\Vert  \label{3.6-}
\end{eqnarray}%
Considering cosine theorem, we have 
\begin{equation}
\left\Vert \overset{\_}{\xi }_{2}(s+\Delta s)-\overset{\_}{\xi }%
_{2}(s)\right\Vert ^{2}=\left\Vert \overset{\_}{\xi }_{2}(s+\Delta
s)\right\Vert ^{2}+\left\Vert \overset{\_}{\xi }_{2}(s)\right\Vert
^{2}-2\left\Vert \overset{\_}{\xi }_{2}(s+\Delta s)\right\Vert \left\Vert 
\overset{\_}{\xi }_{2}(s)\right\Vert \cos h.  \label{3.6**}
\end{equation}%
Since $\overset{\_}{\xi }_{2}(s)$ and $\overset{\_}{\xi }_{2}(s+\Delta s)$
are unit, (\ref{3.6**}) becomes 
\begin{equation*}
\left\Vert \overset{\_}{\xi }_{2}(s+\Delta s)-\overset{\_}{\xi }%
_{2}(s)\right\Vert ^{2}=2-2\cos h.
\end{equation*}%
By using trigonometric relation $\cos h=1-2\sin ^{2}\frac{h}{2},$ last
equality gives $\left\Vert \overset{\_}{\xi }_{2}(s+\Delta s)-\overset{\_}{%
\xi }_{2}(s)\right\Vert =2\sin \frac{h}{2}$. Hence, 
\begin{eqnarray*}
\lim_{\Delta s\rightarrow 0}\left\Vert \frac{\overset{\_}{\xi }_{2}(s+\Delta
s)-\overset{\_}{\xi }_{2}(s)}{\varphi _{1}(s+\Delta s)-\varphi _{1}(s)}%
\right\Vert &=&\lim_{\Delta s\rightarrow 0}\left\vert \frac{1}{\varphi
_{1}(s+\Delta s)-\varphi _{1}(s)}\right\vert \left\Vert \overset{\_}{\xi }%
_{2}(s+\Delta s)-\overset{\_}{\xi }_{2}(s)\right\Vert \\
&=&\lim_{h\rightarrow 0}\frac{1}{h}\left( 2\sin \frac{h}{2}\right) \\
&=&\lim_{h\rightarrow 0}\frac{\sin (h/2)}{h/2} \\
&=&1
\end{eqnarray*}%
Then, from (\ref{3.6-}), we have $p=\left\vert \frac{d\varphi _{1}}{ds}%
\right\vert $. The proof of second equality in (\ref{3.6*}) can be given by
similar way.
\end{proof}

\begin{theorem}
Let $(C,\overset{\_}{\xi })$ be a versor field with curvatures $K_{1},K_{2}$
and alternative curvatures $p,q$. The relations between these curvatures are
given as follows%
\begin{equation}
K_{1}(s)=p(s)\cos \left( \dint q(s)ds\right) ,~\ K_{2}(s)=p(s)\sin \left(
\dint q(s)ds\right) ,  \label{3.6.1}
\end{equation}

\begin{equation}
\begin{array}{c}
a_{1}(s)=-b_{2}(s)\cos \left( \dint q(s)\right) +b_{3}(s)\sin \left( \dint
q(s)ds\right) , \\ 
a_{2}(s)=b_{1}(s), \\ 
a_{3}(s)=b_{2}(s)\sin \left( \dint q(s)ds\right) +b_{3}(s)\cos \left( \dint
q(s)ds\right) .%
\end{array}
\label{3.6.2}
\end{equation}
\end{theorem}

\begin{proof}
Let define the function $\omega (s)=\frac{K_{2}}{K_{1}}(s).$ Then,
alternative curvature $q$ takes the form $q=\frac{\omega ^{\prime }}{%
1+\omega ^{2}}$. Integrating last equality gives $\dint qds=\arctan (\omega
) $. Hence, we obtain $\omega =\tan \left( \dint qds\right) $. Considering
the relation $p=\sqrt{K_{1}^{2}+K_{2}^{2}}$ and trigonometric relations, we
obtained equalities (\ref{3.6.1}). Now, comparing (\ref{2.1}) and (\ref{3.5}%
) and using (\ref{3.6.1}), we have equalities (\ref{3.6.2})
\end{proof}

\begin{definition}
Let $(C,\overset{\_}{\xi })$ be a versor field with AF-frame $R_{A}\{P(s);%
\overset{\_}{\xi }_{2}(s),\overset{\_}{Y}(s),\overset{\_}{D}(s)\}$ in Myller
configuration $M$. The curve $C$ is called Darboux helix in $M$ if the
versor field $\overset{\_}{D}$ makes a constant angle with a constant versor
field $\overset{\_}{l\text{.}}$
\end{definition}

\begin{theorem}
The curve $C$ with AF-frame $R_{A}$ and $\left( K_{1},K_{2}\right) \neq
(0,0) $ in $M$ is Darboux helix iff the function $f$ given below is constant 
\begin{equation}
f=\mp \frac{p}{q}=\mp \frac{1}{\sigma }.  \label{3.7}
\end{equation}
\end{theorem}

\begin{proof}
Let $C$ be a Darboux helix in $M$. Then, there exists a constant angle $\phi 
$ such that $\left\langle \overset{\_}{D},\overset{\_}{l}\right\rangle =\cos
\phi .$ Differentiating last equality with respect to $s$ and considering (%
\ref{3.6}), we obtain $q\left\langle \overset{\_}{Y},\overset{\_}{l}%
\right\rangle =0$. Hence, we have $\overset{\_}{l}\in sp\{\overset{\_}{\xi }%
_{2},\overset{\_}{D}\}$. Now, we can write $\overset{\_}{l}=\mp \sin \phi 
\overset{\_}{\xi }_{2}+\cos \phi \overset{\_}{D}$. Since $\overset{\_}{l}$
is a constant versor field, the differentiation of last equality gives 
\begin{equation*}
(\mp p\sin \phi -q\cos \phi )\overset{\_}{Y}=0.
\end{equation*}%
Hence, we have $\cot \theta =\mp \frac{p}{q}=\mp \frac{1}{\sigma }$ is
constant.
\end{proof}

From Theorem 10, the following corollaries are obtained.

\begin{corollary}
The curve $C$ in $M$ is Darboux helix iff $C$ is slant helix.
\end{corollary}

\begin{corollary}
The axis $\overset{\_}{l}$ of a Darboux helix $C$ in terms of Frenet frame $%
R_{F}$ is given by%
\begin{equation*}
\overset{\_}{l}=\cos \phi \frac{K_{2}}{\sqrt{K_{1}^{2}+K_{2}^{2}}}\overset{\_%
}{\xi }_{1}\mp \sin \phi \overset{\_}{\xi }_{2}+\cos \phi \frac{K_{1}}{\sqrt{%
K_{1}^{2}+K_{2}^{2}}}\overset{\_}{\xi }_{3}.
\end{equation*}
\end{corollary}

\section{Differential equations characterizing slant helices and Darboux
helices in $M$}

In this section, we give general differential equations of a versor field $%
(C,\overset{\_}{\xi })$ in $M$ with respect to both AF-frame $R_{A\text{ }}$
and Frenet-type frame $R_{F}$. Furthermore, we obtain differential equations
characterizing helices, slant helices and Darboux helices in $M$.

\begin{theorem}
Let $(C,\overset{\_}{\xi })$ be a versor field with AF-frame $R_{A}$ and
non-zero alternative curvatures $p,q$. The versor field $\overset{\_}{\xi }%
_{2}$ satisfies the following differential equation 
\begin{equation}
\overset{\_}{\xi }_{2}^{\prime \prime \prime }-\left[ \frac{p^{\prime }}{p}+%
\frac{(pq)^{\prime }}{pq}\right] \overset{\_}{\xi }_{2}^{\prime \prime
}+\left\{ pq\left[ \left( \frac{1}{p}\right) ^{\prime }\frac{1}{q}\right]
^{\prime }+p^{2}+q^{2}\right\} \overset{\_}{\xi }_{2}^{\prime }+pq\left( 
\frac{p}{q}\right) ^{\prime }\overset{\_}{\xi }_{2}=0.  \label{3.7.1}
\end{equation}
\end{theorem}

\begin{proof}
From the second equation in (\ref{3.6}), it follows 
\begin{equation}
\overset{\_}{D}=\frac{1}{q}\overset{\_}{Y}^{\prime }+\frac{p}{q}\overset{\_}{%
\xi }_{2},  \label{3.8}
\end{equation}%
and from the first equation in (\ref{3.6}) we get $\overset{\_}{Y}=\frac{1}{p%
}\overset{\_}{\xi }_{2}^{\prime }$. Considering last equality and third
equation in (\ref{3.6}) we obtain 
\begin{equation}
\overset{\_}{D}^{\prime }=-\frac{q}{p}\overset{\_}{\xi }_{2}^{\prime }.
\label{3.9}
\end{equation}%
Differentiating $\overset{\_}{Y}=\frac{1}{p}\overset{\_}{\xi }_{2}^{\prime }$%
, we have 
\begin{equation}
\overset{\_}{Y}^{\prime }=\left( \frac{1}{p}\right) ^{\prime }\overset{\_}{%
\xi }_{2}^{\prime }+\frac{1}{p}\overset{\_}{\xi }_{2}^{\prime \prime }
\label{3.10}
\end{equation}%
Writing (\ref{3.10}) in (\ref{3.8}) gives 
\begin{equation*}
\overset{\_}{D}=\frac{1}{q}\left( \frac{1}{p}\right) ^{\prime }\overset{\_}{%
\xi }_{2}^{\prime }+\frac{1}{q}\frac{1}{p}\overset{\_}{\xi }_{2}^{\prime
\prime }+\frac{p}{q}\overset{\_}{\xi }_{2}.
\end{equation*}%
Differentiating last equality and by taking into account (\ref{3.9}), we
have (\ref{3.7.1}).
\end{proof}

\begin{corollary}
The curve $C$ with AF-frame $R_{A}$ and non-zero curvatures $p,q$ is slant
helix (or Darboux helix) iff the versor field $\overset{\_}{\xi }_{2}$%
satisfies the following differential equation 
\begin{equation}
\overset{\_}{\xi }_{2}^{\prime \prime \prime }-\left[ \frac{p^{\prime }}{p}+%
\frac{(pq)^{\prime }}{pq}\right] \overset{\_}{\xi }_{2}^{\prime \prime
}+\left\{ pq\left[ \left( \frac{1}{p}\right) ^{\prime }\frac{1}{q}\right]
^{\prime }+p^{2}+q^{2}\right\} \overset{\_}{\xi }_{2}^{\prime }=0.
\label{3.11}
\end{equation}
\end{corollary}

\begin{proof}
From Theorem 10, we have that the curve $C$ is slant helix (or Darboux
helix) iff $\frac{p}{q}$ is constant. Using that in (\ref{3.7.1}), (\ref%
{3.11}) is obtained.
\end{proof}

\begin{theorem}
Let $(C,\overset{\_}{\xi })$ be a versor field with AF-frame $R_{A}$ and
non-zero alternative curvatures $p,q$. The versor field $\overset{\_}{Y}$
satisfies the following differential equation 
\begin{equation*}
\overset{\_}{Y}^{\prime \prime \prime }+\frac{1}{\lambda _{1}}\left( \lambda
_{1}^{\prime }-\lambda _{2}\right) \overset{\_}{Y}^{\prime \prime }+\frac{1}{%
\lambda _{1}}\left( \lambda _{3}-\lambda _{2}^{\prime }\right) \overset{\_}{Y%
}^{\prime }+\frac{1}{\lambda _{1}}\left( \lambda _{3}^{\prime }-p\right) 
\overset{\_}{Y}=0,
\end{equation*}%
where $\lambda _{1}=\frac{q}{q^{\prime }p-qp^{\prime }},$ $\lambda _{2}=%
\frac{1}{p}\left( 1+p^{\prime }\lambda _{1}\right) $ and $\lambda
_{3}=\left( p^{2}+q^{2}\right) \lambda _{1}.$
\end{theorem}

\begin{proof}
From the second equality in (\ref{3.6}) we have 
\begin{equation}
\overset{\_}{\xi }_{2}=-\frac{1}{p}\overset{\_}{Y}^{\prime }+\frac{q}{p}%
\overset{\_}{D}.  \label{3.11c}
\end{equation}%
Differentiating (\ref{3.11c}) and considering first and third equalities in (%
\ref{3.6}), it follows 
\begin{equation}
\overset{\_}{D}=\frac{p^{2}}{q^{\prime }p-p^{\prime }q}\left[ \frac{1}{p}%
\overset{\_}{Y}^{\prime \prime }+\left( \frac{1}{p}\right) ^{\prime }\overset%
{\_}{Y}^{\prime }+\frac{p^{2}+q^{2}}{p}\overset{\_}{Y}\right] .
\label{3.11d}
\end{equation}%
Writing (\ref{3.11d}) in (\ref{3.11c}), we have%
\begin{equation}
\overset{\_}{\xi }_{2}=\frac{q}{q^{\prime }p-p^{\prime }q}\overset{\_}{Y}%
^{\prime \prime }+\left[ \frac{pq}{q^{\prime }p-p^{\prime }q}\left( \frac{1}{%
p}\right) ^{\prime }-\frac{1}{p}\right] \overset{\_}{Y}^{\prime }+\frac{%
q(p^{2}+q^{2})}{q^{\prime }p-p^{\prime }q}\overset{\_}{Y}.  \label{3.11e}
\end{equation}%
Now, differentiating (\ref{3.11e}), we obtain the desired equation.
\end{proof}

\begin{corollary}
The curve $C$ with AF-frame $R_{A}$ and non-zero curvatures $p,q$ is slant
helix (or Darboux helix) iff the versor field $\overset{\_}{Y}$ satisfies
the following differential equation 
\begin{equation*}
\overset{\_}{Y}^{\prime \prime }+q\left( \frac{1}{q}\right) ^{\prime }%
\overset{\_}{Y}^{\prime }+(p^{2}+q^{2})\overset{\_}{Y}=0.
\end{equation*}
\end{corollary}

\begin{proof}
From the second equality in (\ref{3.6}), it follows $\overset{\_}{D}=\frac{1%
}{q}\overset{\_}{Y}^{\prime }+\frac{p}{q}\overset{\_}{\xi }_{2}$ and
differentiating that gives 
\begin{equation*}
\overset{\_}{D^{\prime }}=\left( \frac{1}{q}\right) ^{\prime }\overset{\_}{Y}%
^{\prime }+\frac{1}{q}\overset{\_}{Y}^{\prime \prime }+\left( \frac{p}{q}%
\right) ^{\prime }\overset{\_}{\xi }_{2}+\frac{p}{q}\overset{\_}{\xi }%
_{2}^{\prime }.
\end{equation*}%
By taking into account the first and third equations in (\ref{3.6}), last
equality becomes 
\begin{equation}
\frac{1}{q}\overset{\_}{Y}^{\prime \prime }+\left( \frac{1}{q}\right)
^{\prime }\overset{\_}{Y}^{\prime }+(\frac{p^{2}+q^{2}}{q})\overset{\_}{Y}%
+\left( \frac{p}{q}\right) ^{\prime }\overset{\_}{\xi }_{2}=0.  \label{3.11a}
\end{equation}%
From (\ref{3.11a}), $C$ is a slant helix (or Darboux helix) iff $\overset{\_}%
{Y}^{\prime \prime }+q\left( \frac{1}{q}\right) ^{\prime }\overset{\_}{Y}%
^{\prime }+(p^{2}+q^{2})\overset{\_}{Y}=0$ holds.
\end{proof}

\begin{theorem}
Let $(C,\overset{\_}{\xi })$ be a versor field with AF-frame $R_{A}$ and
non-zero alternative curvatures $p,q$. The versor field $\overset{\_}{D}$
satisfies the following differential equation 
\begin{equation}
\overset{\_}{D}^{\prime \prime \prime }-\left[ \frac{(pq)^{\prime }}{pq}+%
\frac{q^{\prime }}{q}\right] \overset{\_}{D}^{\prime \prime }+\left\{ pq%
\left[ \frac{1}{p}\left( \frac{1}{q}\right) ^{\prime }\right] ^{\prime
}+p^{2}+q^{2}\right\} \overset{\_}{D}^{\prime }+pq\left( \frac{q}{p}\right)
^{\prime }\overset{\_}{D}=0.  \label{3.11.1}
\end{equation}
\end{theorem}

\begin{proof}
From the third equation in (\textbf{\ref{3.6}}), we get $\overset{\_}{Y}=-%
\frac{1}{q}\overset{\_}{D^{\prime }}$ and differentiating that gives 
\begin{equation*}
\overset{\_}{Y^{\prime }}=\left( -\frac{1}{q}\right) ^{\prime }\overset{\_}{%
D^{\prime }}+\left( -\frac{1}{q}\right) \overset{\_}{D^{\prime \prime }.}
\end{equation*}%
Writing that in the second equation in (\ref{3.6}), it follows 
\begin{equation*}
\overset{\_}{\xi }_{2}=\frac{1}{pq}\overset{\_}{D}^{\prime \prime }+\frac{1}{%
p}\left( \frac{1}{q}\right) ^{\prime }\overset{\_}{D}^{\prime }+\frac{q}{p}%
\overset{\_}{D}.
\end{equation*}%
Differentiating last equaliy gives 
\begin{equation}
\overset{\_}{\xi }_{2}^{\prime }=\frac{1}{pq}\overset{\_}{D}^{\prime \prime
\prime }+\left[ \left( \frac{1}{pq}\right) ^{\prime }+\frac{1}{p}\left( 
\frac{1}{q}\right) ^{\prime }\right] \overset{\_}{D}^{\prime \prime }\left\{ %
\left[ \frac{1}{p}\left( \frac{1}{q}\right) ^{\prime }\right] ^{\prime }+%
\frac{q}{p}\right\} \overset{\_}{D}^{\prime }+\left( \frac{q}{p}\right)
^{\prime }\overset{\_}{D}.  \label{3.12}
\end{equation}%
Writing $\overset{\_}{Y}=-\frac{1}{q}\overset{\_}{D^{\prime }}$ in the first
equation in (\ref{3.6}), we have $\overset{\_}{\xi }_{2}^{\prime }=-\frac{p}{%
q}\overset{\_}{D}^{\prime }$. Considering that in (\ref{3.12}), we have (\ref%
{3.11.1}).
\end{proof}

\begin{corollary}
The curve $C$ with AF-frame $R_{A}$ and non-zero curvatures $p,q$ is slant
helix (or Darboux helix) iff the versor field $\overset{\_}{D}$ satisfies
the following differential equation%
\begin{equation}
\overset{\_}{D}^{\prime \prime \prime }-\left[ \frac{(pq)^{\prime }}{pq}+%
\frac{q^{\prime }}{q}\right] \overset{\_}{D}^{\prime \prime }+\left\{ pq%
\left[ \frac{1}{p}\left( \frac{1}{q}\right) ^{\prime }\right] ^{\prime
}+p^{2}+q^{2}\right\} \overset{\_}{D}^{\prime }=0.  \label{3.13}
\end{equation}
\end{corollary}

\begin{proof}
The proof is clear from Theorem 10.
\end{proof}

In Theorem 13, we give the differential equation of a curve $C$ in $M$ with
respect to versor field $\overset{\_}{\xi }_{2}.$ In the following theorems,
we give differential equations of a curve $C$ in $M$ with respect to versor
fields $\overset{\_}{\xi }_{1}$ and $\overset{\_}{\xi }_{3}.$ The proofs of
these theorems can be given by the similar way as given above.

\begin{theorem}
The versor field $\overset{\_}{\xi }_{1}$ of the versor field $(C,\overset{\_%
}{\xi })$ with Frenet frame $R_{F}$ and non-zero curvatures $K_{1},K_{2}$
satisfies the following differential equation%
\begin{equation}
\begin{array}{c}
\overset{\_}{\xi }_{1}^{\prime \prime \prime }-\left[ \frac{\left(
K_{1}K_{2}\right) ^{\prime }}{K_{1}K_{2}}+\frac{K_{1}^{\prime }}{K_{1}}%
\right] \overset{\_}{\xi }_{1}^{\prime \prime } \\ 
+\left\{ K_{1}K_{2}\left[ \frac{1}{K_{2}}\left( \frac{1}{K_{1}}\right)
^{\prime }\right] ^{\prime }+K_{1}^{2}+K_{2}^{2}\right\} \overset{\_}{\xi }%
_{1}^{\prime }+K_{1}K_{2}\left( \frac{K_{1}}{K_{2}}\right) ^{\prime }\overset%
{\_}{\xi }_{1}=0.%
\end{array}
\label{4.1}
\end{equation}%
\qquad\ 
\end{theorem}

\begin{corollary}
The curve $C$ in $M$ with Frenet frame $R_{F}$ and non-zero curvatures $%
K_{1},K_{2}$ is $\overset{\_}{\xi }_{1}$-helix (or $\overset{\_}{\xi }_{3}$%
-helix) iff the following differential equation holds%
\begin{equation}
\overset{\_}{\xi }_{1}^{\prime \prime \prime }-\left[ \frac{\left(
K_{1}K_{2}\right) ^{\prime }}{K_{1}K_{2}}+\frac{K_{1}^{\prime }}{K_{1}}%
\right] \overset{\_}{\xi }_{1}^{\prime \prime }+\left\{ K_{1}K_{2}\left[ 
\frac{1}{K_{2}}\left( \frac{1}{K_{1}}\right) ^{\prime }\right] ^{\prime
}+K_{1}^{2}+K_{2}^{2}\right\} \overset{\_}{\xi }_{1}^{\prime }=0.
\label{4.2}
\end{equation}
\end{corollary}

\begin{theorem}
The versor field $\overset{\_}{\xi }_{2}$of the versor field $(C,\overset{\_}%
{\xi })$  with Frenet frame $R_{F}$ and non-zero curvatures $K_{1},K_{2}$
satisfies the following differential equation%
\begin{equation}
\overset{\_}{\xi }_{2}^{\prime \prime \prime }+\frac{1}{\rho _{1}}\left(
\rho _{1}^{\prime }-\rho _{2}\right) \overset{\_}{\xi }_{2}^{\prime \prime }+%
\frac{1}{\rho _{1}}\left( \rho _{3}-\rho _{2}^{\prime }\right) \overset{\_}{%
\xi }_{2}^{\prime }+\frac{1}{\rho _{1}}\left( \rho _{3}^{\prime
}-K_{1}\right) \overset{\_}{\xi }_{2}=0,  \label{4.2.1}
\end{equation}%
where $\rho _{1}=\frac{K_{2}}{K_{2}^{\prime }K_{1}-K_{2}K_{1}^{\prime }},$ $%
\rho _{2}=\frac{1}{K_{1}}\left( 1+K_{1}^{\prime }\rho _{1}\right) $ and $%
\rho _{3}=\left( K_{1}^{2}+K_{2}^{2}\right) \rho _{1}.$
\end{theorem}

\begin{corollary}
The curve $C$ in $M$ with Frenet frame $R_{F}$ and non-zero curvatures $%
K_{1},K_{2}$ is $\overset{\_}{\xi }_{1}$-helix (or $\overset{\_}{\xi }_{3}$%
-helix) iff the following differential equation holds%
\begin{equation*}
\overset{\_}{\xi }_{2}^{\prime \prime }+K_{1}\left( \frac{1}{K_{1}}\right)
^{\prime }\overset{\_}{\xi }_{2}^{\prime }+(K_{1}^{2}+K_{2}^{2})\overset{\_}{%
\xi }_{2}=0.
\end{equation*}
\end{corollary}

\begin{theorem}
The versor field $\overset{\_}{\xi }_{3}$of the versor field $(C,\overset{\_}%
{\xi })$  with Frenet frame $R_{F}$ and non-zero curvatures $K_{1},K_{2}$
satisfies the following differential equation%
\begin{equation}
\begin{array}{c}
\overset{\_}{\xi }_{3}^{\prime \prime \prime }-\left[ \frac{\left(
K_{1}K_{2}\right) ^{\prime }}{K_{1}K_{2}}+\frac{K_{2}^{\prime }}{K_{2}}%
\right] \overset{\_}{\xi }_{3}^{\prime \prime } \\ 
+\left\{ K_{1}K_{2}\left[ \frac{1}{K_{1}}\left( \frac{1}{K_{2}}\right)
^{\prime }\right] ^{\prime }+K_{1}^{2}+K_{2}^{2}\right\} \overset{\_}{\xi }%
_{3}^{\prime }+K_{1}K_{2}\left( \frac{K_{2}}{K_{1}}\right) ^{\prime }\overset%
{\_}{\xi }_{3}=0.%
\end{array}
\label{4.4}
\end{equation}
\end{theorem}

\begin{corollary}
The curve $C$ in $M$ with Frenet frame $R_{F}$ and non-zero curvatures $%
K_{1},K_{2}$ is $\overset{\_}{\xi }_{1}$-helix (or $\overset{\_}{\xi }_{3}$%
-helix) iff the following differential equation holds%
\begin{equation*}
\overset{\_}{\xi }_{3}^{\prime \prime \prime }-\left[ \frac{\left(
K_{1}K_{2}\right) ^{\prime }}{K_{1}K_{2}}+\frac{K_{2}^{\prime }}{K_{2}}%
\right] \overset{\_}{\xi }_{3}^{\prime \prime }+\left\{ K_{1}K_{2}\left[ 
\frac{1}{K_{1}}\left( \frac{1}{K_{2}}\right) ^{\prime }\right] ^{\prime
}+K_{1}^{2}+K_{2}^{2}\right\} \overset{\_}{\xi }_{3}^{\prime }=0.
\end{equation*}
\end{corollary}

\section{Conclusions}

Slant helices and Darboux helices in Myller configuration $M$ are introduced
and studied. It is shown that a curve $C$ in $M$ is slant helix if and only
if it is Darboux helix in $M$. Moreover, alternative frame of a curve $C$ in 
$M$ is introduced and geometric interpretations of alternative curvatures
are given. Later, general differential equations for a curve $C$ in $M$
according to both alternative frame $R_{A}$ and Frenet type frame $R_{F}$
are obtained. Finally, differential equations characterizing helices, slant
helices and Darboux helices in $M$ \ are given.

\textbf{Author Declaration: }There are no known conflicts of interest
associated with this publication and there has been no significant financial
support for this work that could have influenced its outcome.

\textbf{Data Availability Statement:} This manuscript has no associated data.

\end{document}